\documentclass[12pt]{article}
\usepackage{geometry,graphicx,latexsym,amssymb,verbatim,cite,algorithm2e}
\usepackage{lineno,setspace}
\newcounter{thm}
\newtheorem{theorem}[thm]{Theorem}

\begin{document}


\date{}

\begin{center}
{ \Large On Selkow's  Bound on the Independence Number of Graphs}\\[3mm]
Jochen Harant, Samuel Mohr\\
Ilmenau University of Technology,
Department of Mathematics, Germany 
\end{center}

\begin{abstract}
\noindent
For a graph $G$ with vertex set $V(G)$ and independence number $\alpha(G)$, S. M. Selkow \\(Discrete Mathematics, 132(1994)363--365) established the famous lower bound \\$\sum\limits_{v\in V(G)}\frac{1}{d(v)+1}(1+\max\{\frac{d(v)}{d(v)+1}-\sum\limits_{u\in N(v)}\frac{1}{d(u)+1},0 \})$ on $\alpha(G)$, where $N(v)$ and \\$d(v)=|N(v)|$ denote the neighborhood  and  the degree of a vertex $v\in V(G)$, respectively. However, Selkow's original proof of this result is incorrect. We give a new  probabilistic proof of Selkow's  bound here.
\end{abstract}
We consider a {\it finite}, {\it simple}, and {\it undirected} graph
$G$ with {\it vertex set} $V(G)$.  Let $N_G(v)$ and $d_G(v)=|N_G(v)|$ denote the {\it neighborhood} and the {\it
degree} of $v\in V(G)$, respectively.  A set of vertices $I\subseteq V(G)$  is {\it independent} if no two vertices in $I$ are adjacent.
The {\it independence number} $\alpha(G)$ of $G$ is the maximum
cardinality of an independent set
of $G$.\\
The independence number is one of the most fundamental and
well-studied graph parameters. In view of its
computational hardness,  various bounds on the independence
number have been proposed. The classical lower bound $CW(G)=\sum\limits_{v\in V(G)}\frac{1}{d_G(v)+1}$ on $\alpha(G)$ is due to Y. Caro \cite{Caro} and V. K. Wei \cite{Wei}.
It is natural to ask whether improvements of $\alpha(G)\ge CW(G)$ are possible if more information about $G$ is known than just its degrees.
The following result takes not only the degree of every vertex but also the degree distribution in its neighborhood into account.
\begin{theorem}[S. M. Selkow, \cite{Selkow}]\label{S}$~$\\
$\alpha(G)\ge CW(G)+\sum\limits_{v\in V(G)}\frac{1}{d_G(v)+1}\max\{\frac{d_G(v)}{d_G(v)+1}-\sum\limits_{u\in N_G(v)}\frac{1}{d_G(u)+1},0 \}$.
\end{theorem}
Unfortunately, Selkow's original proof of Theorem \ref{S} is not correct. To our best knowledge, this has not been discovered earlier, and we are not aware of an alternative, correct proof. In order to extract the problematic part of Selkow's argument, let us repeat his proof:\\
 For an event $A$ and a random variable $X$ let $P(A)$ and $E(X)$ denote  be the probability of $A$ and the expectation of $X$, respectively.\\
First, a uniformly chosen ordering $<$ of $V(G)$ is considered.\\
Obviously, the set $I_1=\{v\in V(G)~|~u\in N_G(v)\Rightarrow v<u \}$ is independent and it is easy to show that  $E(|I_1|)=CW(G)$ (e.\,g. see \cite{AS}). \\
Next, let the graph $H$ (depending on the ordering $<$) be obtained from $G$ by removing  $I_1\cup \bigcup\limits_{x\in I_1}N_G(x)$
 and consider the set $I_2=\{v\in V(H)~|~u\in N_H(v)\Rightarrow v<u \}$. \\
 It follows that
$\alpha(G)\ge E(|I_1|+|I_2|)= E(|I_1|)+E(|I_2|)=CW(G)+\sum\limits_{v\in V(G)}P(v\in I_2)$, since $I_1\cap I_2=\emptyset $ and $I_1\cup I_2$ is an independent set of $G$.
To finish the proof of Theorem \ref{S} in \cite{Selkow}, the inequality $P(v\in I_2)\ge \frac{1}{d_G(v)+1}\max\{\frac{d_G(v)}{d_G(v)+1}-\sum\limits_{u\in N_G(v)}\frac{1}{d_G(u)+1},0 \}$ for all $v\in V(G)$  (\cite{Selkow}, page~364, lines 14--17) is used. This turns out to be false in the following general sense.

\noindent{\em For every $\varepsilon>0$, there is a graph $G$ and a vertex $v\in V(G)$, such that \\$0<P(v\in I_2)<\varepsilon \cdot \frac{1}{d_G(v)+1}(\frac{d_G(v)}{d_G(v)+1}-\sum\limits_{u\in N_G(v)}\frac{1}{d_G(u)+1})$.
}\\
To see this, let $n$ be a large positive integer. Consider an arbitrary graph $F$ on $n-3$ vertices and let the graph $G$ on $n$ vertices be obtained by adding three new vertices $v,w,x$ and the edges $vw$, $wx$, and $xy$ for all $y\in V(F)$. For an arbitrary ordering $<$ of $V(G)$,   $v\in I_2$ if and only if $x\in I_1$, $x<w<v$, and $x<y$ for all $y\in V(F)$.
 It is easy to see that there are exactly ${n-1 \choose 2}(n-3)!$ such orderings with the property $v\in I_2$, thus, $P(v\in I_2)=\frac{{n-1 \choose 2}(n-3)!}{n!}=\frac{1}{2n}$, however, $\frac{1}{d_G(v)+1}(\frac{d_G(v)}{d_G(v)+1}-\sum\limits_{u\in N_G(v)}\frac{1}{d_G(u)+1})=\frac12(\frac12-\frac13)=\frac{1}{12}$. \\ \\
Eventually, we present a new probabilistic proof of Theorem \ref{S}.  \\ \\
{\bf Proof of Theorem \ref{S}.}
As in Selkow's proof, consider a uniformly chosen ordering $<$ of $V(G)$, the set $I_1=\{v\in V(G)~|~u\in N_G(v)\Rightarrow v<u \}$, and the  graph $H$ induced by $V(G)\setminus (I_1\cup \bigcup\limits_{x\in I_1}N_G(x))$.
\\ With  $f(v)=0$ if $v\notin V(H)$  and $f(v)=\frac{1}{d_G(v)+1}$ if $v\in V(H)$, it follows \\
$\alpha(G)\ge E(|I_1|+\alpha(H)) = E(|I_1|)+E(\alpha(H)) \ge E(|I_1|)+E(CW(H))= CW(G)+E(\sum\limits_{v\in V(H)}\frac{1}{d_H(v)+1})$\\
$\ge CW(G)+E(\sum\limits_{v\in V(H)}\frac{1}{d_G(v)+1})=CW(G)+E(\sum\limits_{v\in V(G)}f(v))$\\
$=CW(G)+\sum\limits_{v\in V(G)}\frac{1}{d_G(v)+1}P(v\in V(H))$. \\
Using
$P(v\in V(H))\ge 0$ and
  $P(v\notin V(H))=P(v\in I_1\vee (\bigvee\limits_{u\in N(v)}u\in I_1))$\\
  $\le P(v\in I_1)+\sum\limits_{u\in N(v)}P(u\in I_1) =\frac{1}{d_G(v)+1}+\sum\limits_{u\in N(v)}\frac{1}{d_G(u)+1}$,
    Theorem \ref{S} is proved. \hfill $\Box$

\end{document}